\RequirePackage{fix-cm}
\documentclass[smallextended]{svjour3}       
\smartqed  
\usepackage{graphicx}
\begin{document}

\title{A new Improvement of Ambarzumyan's Theorem
}

\titlerunning{An inverse result}        

\author{Alp Arslan K\i ra\c{c}
}


\institute{A. A. K\i ra\c{c} \at
               Department of Mathematics, \\ Faculty of Arts and Sciences, Pamukkale
University, \newline
 20070, Denizli, Turkey\\
              \email{aakirac@pau.edu.tr}           
}

\date{Received: date / Accepted: date}

\maketitle

\begin{abstract} We extend the classical Ambarzumyan's theorem to the quasi-periodic boundary value problems by using only a part knowledge of one
spectrum. We also weaken slightly the Yurko's conditions on the first eigenvalue.       
\end{abstract}

\section{Preliminary Improvement}
In $1929$, Ambarzumyan \cite{Ambarz} proved that if  $\{(n\pi)^{2}: n=0,1,2\ldots\}$ is the spectrum of the boundary value problem
 \begin{equation}  \label{1}
-y^{\prime\prime}(x)+q(x)y(x)=\lambda y(x),\qquad y^{\prime}(0)=y^{\prime}(1)=0
\end{equation}
with real potential $q\in L^{1}[0,1]$, then $q=0$ a.e. Clearly, if $q=0$ a.e., then the eigenvalues $\lambda_{n}=(n\pi)^{2}$, $n\geq 0$.

Freiling and Yurko \cite{FreilingYurko} proved that it is enough to specify only the first eigenvalue rather than the whole spectrum. More precisely, the first eigenvalue denoted by $\lambda_{0}$ is a mean value of the
potential, that is, they proved the following Ambarzumyan-type theorem:

\begin{theorem} \label{yurk1}
If $\lambda_{0}=\int_{0}^{1}q(x)\,dx$, then $q=\lambda_{0}$ a.e.
\end{theorem}

In \cite{Yurko2013}, Yurko also provided generalizations of Theorem \ref{yurk1} on wide classes of self-adjoint differential operators. Some of the inverse results in \cite{Yurko2013} are in the following theorem:

\begin{theorem} \label{yurk2}
(a) Let
\[
\lambda_{0}=\int_{0}^{1}q(x)\,dx
\]
be the first eigenvalue of the periodic boundary value problem
\begin{equation}  \label{per}
-y^{\prime\prime}(x)+q(x)y(x)=\lambda y(x),\qquad y(1)=y(0),\quad y^{\prime}(1)=y^{\prime}(0),
\end{equation}
 then $q=\lambda_{0}$ a.e.

 (b) Let
 \begin{equation}  \label{condper+}
\lambda_{0}=\pi^{2}+\frac{2}{\alpha^2+\beta^2}\int_{0}^{1}q(x)(\alpha sin\pi x +\beta cos\pi x)^2\,dx,
\end{equation}
 for some fixed $\alpha$ and $\beta$, be the first eigenvalue of the anti-periodic boundary value problem
\begin{equation}  \label{per+}
-y^{\prime\prime}(x)+q(x)y(x)=\lambda y(x),\qquad y(1)=-y(0),\quad y^{\prime}(1)=-y^{\prime}(0).
\end{equation}
 Then \[q=\frac{2}{\alpha^2+\beta^2}\int_{0}^{1}q(x)(\alpha sin\pi x +\beta cos\pi x)^2\,dx\,\,\,\, a.e.\]
\end{theorem}

Consider the boundary value problems $L_{t}(q)$ generated in the space $L^{2}[0, 1]$ by the following differential equation
and quasi-periodic boundary conditions
\begin{equation}\label{q-per}
-y^{\prime\prime}(x)+q(x)y(x)=\lambda y(x),\qquad y(1)=e^{it}y(0),\quad y^{\prime}(1)=e^{it}y^{\prime}(0),
\end{equation}
where $q\in L^{1}[0,1]$ is a real-valued function and  $t\in[0,\pi)\cup [\pi,2\pi)$.
The operator $L_{t}(q)$ is  self-adjoint and the cases $t=0$ and $t=\pi$ correspond to the periodic and anti-periodic problems, respectively. Let $\{\lambda_{n}(t)\}_{n\in\mathcal{Z}}$ be the eigenvalues of the operator $L_{t}(q)$. In the case when $q=0$, $(2\pi n+t)^{2}$ for $n\in \mathcal{Z}$ is the eigenvalue of the operator $L_{t}(0)$ for any fixed $t\in [0,2\pi)$ corresponding to the eigenfunction $e^{i(2\pi n+t)x}$.

The result of Ambarzumyan \cite{Ambarz} is an exception to the general rule. In general, Borg \cite{Borg} proved that one spectrum
does not determine the potential. Also, Borg showed that two spectra determine it uniquely. Many generalizations of the Ambarzumyan's theorem can be found in \cite{Chern,corri,Lieberman,Levitan,Ambarzcoupled,anoteinver}. But as far as we know, for the first time, we extend the classical Ambarzumyan's theorem to the quasi-periodic boundary value problems $L_{t}(q)$  by using only a part knowledge of one
spectrum (see Theorem \ref{quasi2}). More precisely,
the aim of this paper is to weaken slightly the conditions in Theorem \ref{yurk1} and Theorem \ref{yurk2} for the first eigenvalue (see Theorem \ref{quasi1}) and to prove a new generalization of Ambarzumyan's theorem, without using any additional conditions on the potential. We also extend Theorem 1 of the previous paper \cite{kýrac:ambars}.

Our new uniqueness-type results read as follows:

\begin{theorem} \label{quasi1} Let $\lambda_{0}(t)$ be the first eigenvalue of $L_{t}(q)$. Then:

(a) If $\lambda_{0}(t)\geq t^{2}+\int_{0}^{1}q(x)\,dx$ for $t\in[0,\pi)$, then $q=\int_{0}^{1}q(x)\,dx$ a.e.

(b) If $\lambda_{0}(t)\geq (2\pi-t)^{2}+\int_{0}^{1}q(x)\,dx$ for $t\in[\pi,2\pi)$, then $q=\int_{0}^{1}q(x)\,dx$ a.e.
\end{theorem}

\begin{theorem} \label{quasi2} Let $\lambda_{0}(t)$ be the first eigenvalue of $L_{t}(q)$, and suppose that $n_{0}$ is a sufficiently large positive integer. Then the following two assertions hold:

(a) If  $\lambda_{0}(t)\geq t^{2}$ and $\lambda_{n}(t)=(2n\pi+t)^{2}$ for $t\in[0,\pi)$ and $n>n_{0}$, then $q=0$ a.e.

(b) If  $\lambda_{0}(t)\geq (2\pi-t)^{2}$ and $\lambda_{n}(t)=(2n\pi+t)^{2}$ for $t\in[\pi,2\pi)$ and $n>n_{0}$, then $q=0$ a.e.
\end{theorem}
In Theorem \ref{quasi1} (a), the case $t=0$ corresponds to the periodic problem (\ref{per}) and the assertion of Theorem \ref{yurk2} (a) holds by using $\lambda_{0}(0)\geq\int_{0}^{1}q(x)\,dx$ instead of $\lambda_{0}(0)=\int_{0}^{1}q(x)\,dx$. In Theorem \ref{quasi1} (b), the case $t=\pi$ corresponds to the anti-periodic problem (\ref{per+}) and the assertion of Theorem \ref{quasi1} (b) holds by using $\lambda_{0}(\pi)\geq\pi^{2}+\int_{0}^{1}q(x)\,dx$ instead of (\ref{condper+}), namely the first eigenvalue depends only on the mean value of the
potential as in Theorem \ref{yurk1} and Theorem \ref{yurk2} (a).

 And, for the boundary value problem (\ref{1}), the form with reduced spectrum of Ambarzumyan's theorem read as follows:
\begin{theorem} \label{neum}
(a) If $\lambda_{0}\geq \int_{0}^{1}q(x)\,dx$, then $q=\int_{0}^{1}q(x)\,dx$ a.e.

(b) If $\lambda_{0}\geq 0$ and $\lambda_{n}=(n\pi)^{2}$ for $n>n_{0}$ , then $q=0$ a.e.,
where  $n_{0}$ is a sufficiently large positive integer.
\end{theorem}
Note that, for example, in Theorem \ref{yurk1} if the first eigenvalue  $\lambda_{0}=0$ and $\int_{0}^{1}q(x)\,dx=0$, then $q=0$ a.e. The first eigenvalue-type of Ambarzumyan's theorem, such as Theorem \ref{yurk1}-\ref{quasi1} and \ref{neum} (a), depends on a mean value of the
potential. Hence, to prove Theorem \ref{quasi2} and  \ref{neum} (b), without imposing
an additional condition on the potential such as $\int_{0}^{1}q(x)\,dx=0$, we have information about the first eigenvalue with a less restrictive one and a subset of the
sufficiently large eigenvalues of the spectrum.

\section{Proofs}

\textbf{Proof of Theorem \ref{quasi1}.} (a) We show that $y=e^{itx}$ is the first eigenfunction corresponding to the first eigenvalue $\lambda_{0}(t)$ of the operator $L_{t}(q)$ for $t\in[0,\pi)$. Since the test function $y=e^{itx}$ satisfies the quasi-periodic boundary conditions in (\ref{q-per}), by the variational principle, we get

\begin{equation}  \label{ineq1}
t^{2}+\int_{0}^{1}q(x)\,dx\leq\lambda_{0}(t)\leq \frac{\int_{0}^{1}-\bar{y}y^{\prime\prime}dx+\int_{0}^{1}q(x)|y|^{2}dx}{\int_{0}^{1}|y|^{2}dx}=t^{2}+\int_{0}^{1}q(x)\,dx.
\end{equation}
This implies that $\lambda_{0}(t)=(t^{2}+\int_{0}^{1}q(x)\,dx)$ is the first eigenvalue corresponding to the first eigenfunction $y=e^{itx}$. Substituting this into the equation
 $$-y^{\prime\prime}+q(x)y=\lambda y,$$
we get $q=\int_{0}^{1}q(x)\,dx$ a.e.

(b) Similarly, for $t\in[\pi,2\pi)$,  the test function $y=e^{i(-2\pi+t)x}$ is the first eigenfunction corresponding to the first eigenvalue $\lambda_{0}(t)=(2\pi-t)^{2}+\int_{0}^{1}q(x)\,dx$. Thus, $q=\int_{0}^{1}q(x)\,dx$ a.e.
$\qed$\\\\
\textbf{Proof of Theorem \ref{quasi2}.} Note that by \cite{Melda.O} (see also \cite{Veliev:Kiraç}), without using the assumption $\int_{0}^{1}q(x)\,dx=0$, the eigenvalues $\lambda_{n}(t)$ of the operator $L_{t}(q)$ for $t\neq 0,\pi$ are asymptotically located  such that
\begin{equation}  \label{melvel}
\lambda_{n}(t)=(2\pi n+t)^{2}+ \int_{0}^{1}q(x)\,dx +O\left(n^{-1}ln|n|\right),
\end{equation}
for $|n|\geq n_{0}$, where $n_{0}$ is a sufficiently large positive integer. And, by Theorem 1.1 of \cite{anoteinver}, there is a similar asymptotic formulas for  the sufficiently large eigenvalues $\lambda_{n}(t)$ of the operator $L_{t}(q)$ for  $t=0,\pi$. Thus, for all $t\in[0,\pi)\cup [\pi,2\pi)$, if $\lambda_{n}(t)=(2n\pi+t)^{2}$ for $n>n_{0}$, then $\int_{0}^{1}q(x)\,dx=0$. From Theorem \ref{quasi1}, $q=0$ a.e. Thus Theorem \ref{quasi2} (a) and (b) are proved. $\qed$
\\\\
\textbf{Remark.} $n<0$ implies that $-n>0$. In Theorem \ref{quasi2} (a) and (b), if we use the large eigenvalues $\lambda_{-n}(t)=(2n\pi-t)^{2}$ for $n>n_{0}$ instead of $\lambda_{n}(t)=(2n\pi+t)^{2}$, then the assertions of Theorem \ref{quasi2} remain valid. Moreover, for all $n>n_{0}$, it is enough to set either of the eigenvalues $\lambda_{n}(t)$,  $\lambda_{-n}(t)$.
\\\\
\textbf{Proof of Theorem \ref{neum}.} (a) Arguing as in the proof of Theorem \ref{quasi1} (a), for $\lambda_{0}\geq \int_{0}^{1}q(x)\,dx$, the test function $y=1$ is the first eigenfunction corresponding to the first eigenvalue $\lambda_{0}=\int_{0}^{1}q(x)\,dx$. Thus, $q=\int_{0}^{1}q(x)\,dx$ a.e.

(b) It follows from the asymptotic formula (1.6.6) in \cite{Levitan} that if
$\lambda_{n}=(n\pi)^{2}$ for $n>n_{0}$, then certainly $\int_{0}^{1}q(x)\,dx=0$. From (a), $q=0$ a.e. \qed


\end{document}